\documentclass[12pt]{amsart}

\usepackage{amsmath, amsthm, amscd, amsfonts, amssymb, graphicx, color}
\usepackage{times}
\usepackage{bm}

\RequirePackage[numbers]{natbib}

\begin{document}

\title{Quantifying nuisance parameter effects via decompositions of asymptotic refinements for likelihood-based statistics}

\author{THOMAS J. DICICCIO}
\address{Department of Social Statistics,
Cornell University, Ithaca, New York 14853, U.S.A.} \email{tjd9@cornell.edu}

\author{TODD A. KUFFNER}
\address{Department of Mathematics, Washington University in St. Louis, St. Louis, Missouri 63130, U.S.A.} \email{kuffner@math.wustl.edu}

\author{G. ALASTAIR YOUNG}
\address{Department of Mathematics, Imperial College London, London SW7 2AZ,
U.K.} \email{alastair.young@imperial.ac.uk}

\begin{abstract}
Accurate inference on a scalar interest parameter in the presence of a nuisance parameter may be obtained using an adjusted version of the signed root likelihood ratio statistic, in particular Barndorff-Nielsen's $R^*$ statistic. The adjustment made by this statistic may be decomposed into a sum of two terms, interpreted as correcting respectively for the possible effect of nuisance parameters and the deviation from standard normality of the signed root likelihood ratio statistic itself. We show that the adjustment terms are determined to second-order in the sample size by their means. Explicit expressions are obtained for the leading terms in asymptotic expansions of these means. These are easily calculated, allowing a simple way of quantifying and interpreting the respective effects of the two adjustments, in particular of the effect of a high dimensional nuisance parameter. Illustrations are given for a number of examples, which provide theoretical insight to the effect of nuisance parameters on parametric inference. The analysis provides a decomposition of the mean of the signed root statistic involving two terms: the first has the property of taking the same value whether there are no nuisance parameters or whether there is an orthogonal nuisance parameter, while the second is zero when there are no nuisance parameters. Similar decompositions are discussed for the Bartlett correction factor of the likelihood ratio statistic, and for other asymptotically standard normal pivots.
\end{abstract}


\maketitle

\section{Introduction}

We are concerned with inference on a scalar interest parameter in the presence of a, possibly high dimensional, nuisance parameter, based on a data sample of size $n$, and with identification of procedures which yield repeated sampling accuracy. In this setting, inference accurate to third order, that is with repeated sampling error of order $O(n^{-3/2})$,  may be obtained using an adjusted version of the signed root likelihood ratio statistic, in particular through use of Barndorff-Nielsen's $R^*$ statistic (Barndorff-Nielsen, 1986).

The $R^*$ statistic is particularly useful in two contexts. In full, multi-parameter exponential family models inference based on standard normal approximation to the sampling distribution of the $R^*$ statistic approximates to third order the optimal, conditional, but generally intractable, inference, which is based on conditioning on the sufficient statistic for the nuisance parameter. In more general models which admit an ancillary statistic, taken to mean an approximately distribution free statistic which together with the maximum likelihood estimator constitutes a minimal sufficient statistic for the full parameter in the model, the normal approximation approximates to the same third order an exact inference based on conditioning on the ancillary statistic. A practical limitation of the use of $R^*$ is in the requirement of explicit specification of the appropriate ancillary, and the need to express the likelihood directly in terms of the maximum likelihood estimator and the ancillary statistic. When calculation of the $R^*$ statistic is tractable, inference with repeated sampling accuracy $O(n^{-3/2})$ is obtained through the normal approximation. This same level of repeated sampling accuracy may be obtained by parametric bootstrap procedures, in particular those based on simulation estimation of the sampling distribution of the unadjusted signed root statistic: see DiCiccio et al.\ (2001), Lee \& Young (2005). Key to this bootstrap approach is appropriate handling of the nuisance parameter: third order repeated sampling accuracy is obtained by considering the sampling distribution of the signed root statistic when the nuisance parameter is specified as the constrained maximum likelihood value calculated from the observed data sample.

Inference based on the $R^*$ statistic and the parametric bootstrap alternative sketched above are analytically related. DiCiccio \& Young (2008) observe that in the problem of inference on a scalar component of the canonical parameter in the multi-parameter exponential family context, inference based on normal approximation to $R^*$ may be viewed as an analytic, saddlepoint approximation to the bootstrap inference. In the same way, it is readily seen that in the ancillary statistic context, inference based on $R^*$ may be regarded as a saddlepoint approximation to a conditional bootstrap calculation, which simulates the distribution of the signed root statistic conditional on the observed value of the ancillary statistic, with the nuisance parameter fixed at its constrained maximum likelihood value. Simulation of this conditional bootstrap distribution will be infeasible in many circumstances, though in certain cases, such as regression-scale models, simple methods of conditional simulation, employing MCMC, are possible: see Brazzale \& Davison (2008). Alternatively, and more simply, the conditional distribution may be replaced by simulation of the marginal distribution of the signed root statistic. DiCiccio et al.\ (2015) demonstrate that the marginal bootstrap distribution approximates the conditional bootstrap distribution to second order, $O(n^{-1})$, given the ancillary statistic.

The adjustment made by the $R^*$ statistic may be decomposed into a sum of two terms, interpreted as correcting respectively for the possible effect of nuisance parameters and an information adjustment, representing the deviation from standard normality of the signed root likelihood ratio statistic itself.  Pierce \& Peters (1992) proposed such a decomposition in the case where the interest parameter is a component of the canonical parameter in a full exponential family model. A generalization of the decomposition is detailed by Barndorff-Nielsen \& Cox (1994, Section 6.6.4). Starting from numerical investigations by Pierce \& Peters (1992), it has been noted that the information adjustment is typically small when the adjusted information for the interest parameter, which we define formally in Section 2, is large. By contrast, the nuisance parameter adjustment can be appreciable when information on the nuisance parameter is small, as will usually occur when its dimension is large. Crucially, however, the magnitude of the nuisance parameter adjustment relative to the information adjustment also depends on the structure of the statistical model in question, and a simple methodology for measurement of nuisance parameter effects for a given model is lacking.

In this paper we note that the adjustment terms are, from a repeated sampling perspective, determined to second-order, $O(n^{-1})$, in the sample size by their means. The precise definitions of the adjustment terms themselves are unimportant to our strategy for quantifying nuisance parameter effects, though we note that, except for full exponential family and transformation models, they must generally be approximated, leading to only second-order accuracy from the resulting adjusted signed root statistic. Approximations to $R^*$ which yield second-order accuracy include those described by DiCiccio \& Martin (1993) and Skovgaard (1996): for a summary see Severini (2000, Section 7.5).

We obtain explicit expressions for the leading terms in asymptotic expansions of the repeated sampling means of the nuisance parameter and information adjustments. These involve calculation only of expectations of certain low-order log-likelihood derivatives, and are therefore easily evaluated for quite general models, even when the $R^*$ statistic itself is intractable. The adjustment terms have variances of low order $O(n^{-2})$ and the asymptotic means therefore allow a simple, effective and general way of quantifying and interpreting the respective effects of the two adjustments. Of particular methodological interest is analysis of the effect of a high dimensional nuisance parameter on the inference based on the $R^*$ statistic, and by extension its bootstrap alternative. Inference based on the $R^*$ statistic, when tractable, represents a `gold standard' in what is achievable in the inference problem and we have noted a close relationship between inference based on the $R^*$ statistic and parametric bootstrap inference. It is reasonable therefore to expect that the calculations are useful too in shedding light on operation of the parametric bootstrap. The repeated sampling properties of the bootstrap are, modulo Monte Carlo error introduced by the need in practice to construct the bootstrap estimate of the sampling distribution of the signed root statistic from a finite simulation, determined entirely by nuisance parameter effects, through substitution of unknown values by estimates. A central recommendation of this paper is that valuable insights to operation of the parametric bootstrap may be obtained by identification of the explicit way in which the means of the nuisance parameter and information adjustments depend on the nuisance parameter. As we shall see in Section 4, in certain key problems these quantities depend only on the dimension of the nuisance parameter, and not on its actual value. In such cases we may reasonably expect good repeated sampling accuracy from the bootstrap, as precise specification of the nuisance parameter values in the calculation is unimportant. In other situations, we observe that the value of the nuisance parameter has a more substantial effect on the adjustment means, in which case we may be alert to impaired accuracy from the bootstrap and its analytic alternatives, especially with small sample sizes.

Our analysis provides a decomposition of the mean of the signed root statistic involving two terms: the first has the property of taking the same value whether there are no nuisance parameters or whether there is an orthogonal nuisance parameter, while the second is zero when there are no nuisance parameters. Similar decompositions are discussed for the Bartlett correction factor of the likelihood ratio statistic, and for other asymptotically standard normal pivots, in Sections 5 and 6 respectively.

\section{The inferential problem}

Suppose that $Y=(Y_1,\ldots,Y_n)$ is a continuous random vector
and that the distribution of $Y$ depends on an unknown $d$-dimensional parameter $\theta=(\theta^1, \ldots, \theta^d)$,
partitioned as $\theta=(\psi,\phi)$, where $\psi=\theta^1$ is a scalar interest parameter
and $\phi$ is a nuisance parameter of dimension $d-1$.
Let $L(\theta)$ be the log-likelihood function for $\theta$ based on $Y$
and let $\hat\theta=(\hat\psi,\hat\phi)$ be the global maximum likelihood estimator of $\theta$.
Further, let $\tilde\theta=\tilde\theta(\psi)=(\psi,\tilde\phi)=\{\psi,\tilde\phi(\psi)\}$ be the constrained maximum likelihood estimator of $\theta$ for given $\psi$. Then the profile log-likelihood function for $\psi$ is $M(\psi)=L\{\tilde\theta(\psi)\}$
and the likelihood ratio statistic for $\psi$ is $W(\psi)=2\{M(\hat\psi)-M(\psi)\}$, where $M(\hat\psi)=L(\hat\theta)$, since $\tilde\theta(\hat\psi)=\hat\theta$.
The signed root likelihood ratio statistic is $R(\psi)={\rm sgn}(\hat\psi-\psi)\{W(\psi)\}^{1/2}$. Then, for example, testing $H_0:\psi=\psi_0$ against $H_a:\psi>\psi_0$ or $H_a:\psi<\psi_0$ can be based on the test statistic $R(\psi_0)$.
Asymptotically, as the sample size $n$ increases, the sampling distribution of $R(\psi)$ tends to the standard normal distribution. Specifically, $R(\psi)$ is distributed as standard normal to first order, to error of order $O(n^{-1/2})$. By contrast, the $R^*$ statistic is distributed as standard normal to error of order $O(n^{-3/2})$.

The $R^*$ statistic is defined by
\begin{equation}
\label{definitionrstar}
R^*(\psi)=R(\psi)+R(\psi)^{-1}\log(v(\psi)/R(\psi)),
\end{equation}
where $v(\psi)$
is given (Barndorff-Nielsen, 1986) by

\begin{equation}
\label{definitioncorrection}
v(\psi)=\left| \begin{array}{cc}
L_{;\hat \theta}(\hat \theta)-L_{;\hat \theta}(\tilde \theta) \\
L_{\phi ;\hat \theta}(\tilde \theta)
\end{array} \right| /\{|j_{\phi \phi}(\tilde \theta)|^{1/2}
|j(\hat \theta)|^{1/2}\}.
\end{equation}
Here, it is supposed that the log-likelihood
function has been written as $L(\theta; \hat \theta, a)$, with $(\hat \theta, a)$
minimal sufficient and $a$ ancillary, that is with a distribution which, at least approximately, does not depend on $\theta$. Further,
\[
L_{; \hat\theta}(\theta)\equiv L_{; \hat \theta}(\theta; \hat \theta, a)=\frac{\partial}{\partial \hat \theta}
L(\theta ;\hat \theta,a), \; \; L_{\phi ;\hat\theta}(\theta)\equiv L_{\phi ;\hat \theta}(\theta;\hat\theta,a)=\frac{\partial^2}
{\partial \phi \partial \hat \theta}L(\theta ;\hat \theta, a).
\]
Also, $j$ denotes the observed information matrix, $j(\theta)=(-L_{rs}(\theta))$, with
$L_{rs}(\theta)=\partial^2 L(\theta)/\partial\theta^r\partial\theta^s$,
and $j_{\phi \phi}$ denotes its $(\phi,\phi)$
component. The sampling distribution of $R^*(\psi)$ is standard normal conditionally on $a$, and hence, as noted,
unconditionally, to error of third order $O(n^{-3/2})$. Note that in a full exponential family model, $\hat \theta$ is already itself sufficient, and no ancillary statistic $a$ is required. The expression for $v(\psi)$ given by (\ref{definitioncorrection}) therefore simplifies somewhat: see, for example, Barndorff-Nielsen \& Cox (1994, Example 6.19).

Barndorff-Nielsen \& Cox (1994, Section 6.6.4), generalizing Pierce \& Peters (1992), introduce quantities ${\rm NP}(\psi)$ and ${\rm INF}(\psi)$, both of order $O_p(n^{-1/2})$, such that $R^{*}(\psi)=R(\psi)+{\rm NP}(\psi)+{\rm INF}(\psi)$. Explicitly, we have
\[
{\rm NP}(\psi)=-\frac{1}{R(\psi)}\log C(\psi),
\]
where
\[
C(\psi)=\frac{\{|j_{\phi \phi}(\hat \theta)||j_{\phi \phi}(\tilde \theta)|\}^{1/2}}{|L_{\phi;\hat \phi}(\tilde \theta)|},
\]
with $L_{\phi ;\hat\phi}(\theta)\equiv L_{\phi ;\hat \phi}(\theta;\hat\theta,a)={\partial^2}
L(\theta ;\hat \theta, a)/\partial \phi \partial \hat \phi$ and, as before, $j_{\phi \phi}$ denoting the $(\phi,\phi)$
component of the observed information $j$.
Also,
\[
{\rm INF}(\psi)=\frac{1}{R(\psi)}\log \{u(\psi)/R(\psi)\},
\]
where
\[
u(\psi)=j_p(\hat \psi)^{-1/2}\frac{\partial}{\partial \hat \psi}\{M(\hat \psi)-M(\psi)\}.
\]
Here $j_p$ is the profile observed information, $j_p(\psi)=-\partial^2 M(\psi)/\partial \psi^2$, and the derivative with respect to $\hat \psi$ is calculated with $M(\hat \psi)-M(\psi)$ considered as a function of $\psi, \hat \psi, \tilde \phi(\psi)$ and $a$.

Calculation of $R^*(\psi)$ supposes explicit representation of the log-likelihood as a function of $(\hat \theta, a)$. Other formulations of the adjustment $v(\psi)$, due to Fraser and co-workers, are possible. The tangent exponential model introduced  by Fraser (1990) avoids the need to specify the transformation $Y \to (\hat \theta, a)$, though still requires awkward analytic calculation: a useful summary is given by Brazzale et al.\ (2007, Chapter 8). In general, however, it is necessary to approximate to the quantity $v(\psi)$. Replacing $v(\psi)$ in the definition (\ref{definitionrstar}) of $R^*(\psi)$ by an estimate $\tilde v(\psi)$ typically yields an adjusted version of the signed root likelihood ratio statistic distributed as standard normal only to error of second order, $O(n^{-1})$. A computationally attractive approximation based on orthogonal parameterisation (Cox \& Reid, 1987) is described by DiCiccio \& Martin (1993). The approximation due to Skovgaard (1996) is theoretically attractive in that it also provides large deviations protection.

To develop our analysis, some further notation is required.  Let $L_\theta(\theta)$ denote the score function, the vector with components $L_r(\theta)=\partial L(\theta)/\partial\theta^r, r=1, \ldots, d$.
In the calculations that follow, arrays and summation are denoted by using the standard conventions, for which the indices $r,s,t,\ldots$ are assumed to range over $1,\ldots,d$. Summation over the range is implied for any index appearing in an expression both as a subscript and as a superscript.
As above, differentiation is indicated by subscripts.
Then $E\{L_r(\theta)\}=0$; let $\lambda_{rs}=E\{L_{rs}(\theta)\}$, $\lambda_{rst}=E\{L_{rst}(\theta)\}$, etc., and
put $l_r=L_r(\theta)$, $l_{rs}=L_{rs}(\theta)-\lambda_{rs}$, $l_{rst}=L_{rst}(\theta)-\lambda_{rst}$, etc.
The constants $\lambda_{rs}$, $\lambda_{rst}, \ldots$, are assumed to be of order $O(n)$.
The variables $l_r$, $l_{rs}$, $l_{rst}$, etc., each of which have expectation 0, are assumed to be of order $O_p(n^{1/2})$.
The joint cumulants of $l_r$, $l_{rs}$, etc.\ are assumed to be of order $O(n)$. These assumptions will usually be satisfied in situations involving independent
observations, or structured dependence, such as in time series contexts.
It is useful to extend the $\lambda$-notation: let $\lambda_{r,s}=E(L_rL_s)=E(l_rl_s)$, $\lambda_{rs,t}=E(L_{rs}L_t)=E(l_{rs}l_t)$, etc.
Bartlett identities involving the $\lambda$'s can be derived by repeated differentiation of the identity $\int\exp\{L(\theta)\}dy=1$; in particular,
$$
\lambda_{rs}+\lambda_{r,s}=0, \quad
\lambda_{rst}+\lambda_{rs,t}+\lambda_{rt,s}+\lambda_{st,r}+\lambda_{r,s,t}=0.
$$
Differentiation of the definition $\lambda_{rs}=\int L_{rs}(\theta)\exp\{L(\theta)\}dy$ yields
$\lambda_{rs/t}=\lambda_{rst}+\lambda_{rs,t},$
where $\lambda_{rs/t}=\partial\lambda_{rs}/\partial\theta^t$.
Further, let $(\lambda^{rs})$ be the $d \times d$ matrix inverse of $(\lambda_{rs})$, and let $\eta=-1/\lambda^{11}$, $\tau^{rs}=\eta\lambda^{1r}\lambda^{1s}$,  and $\nu^{rs}=\lambda^{rs}+\tau^{rs}$.
Thus, $\lambda^{rs}$, $\tau^{rs}$, and $\nu^{rs}$ are of order $O(n^{-1})$, while $\eta$, which is what we have termed the adjusted information for $\psi$, is of order $O(n)$.

DiCiccio \& Stern (1994a) showed that $R(\psi)=\eta^{1/2}\{R_1+R_2+O_p(n^{-3/2})\}$, where $R_1=-\lambda^{1r}l_r$ and
\[
R_2=\lambda^{1r}\lambda^{st}l_{rs}l_t+{\textstyle{1\over
2}}\lambda^{1r}\tau^{st}l_{rs}l_t -{\textstyle{1\over
2}}\lambda^{1r}\lambda^{su}\nu^{tv}\lambda_{rst}l_ul_v
-{\textstyle{1\over 6}}\lambda^{1r}\tau^{su}\tau^{tv}\lambda_{rst}l_ul_v.
\]
Note that $R_1$ is of order $O_p(n^{-1/2})$ and $R_2$ is of order $O_p(n^{-1})$.
Since $E(R_1)=0$, it follows that
\begin{equation}
\label{expectR}
E\{R(\psi)\}=\eta^{1/2}\{\lambda^{1r}\lambda^{st}\lambda_{rs,t}
+{\textstyle{1\over 2}}\lambda^{1r}\tau^{st}\lambda_{rs,t}
+{\textstyle{1\over 2}}\lambda^{1r}\lambda^{st}\lambda_{rst}
+{\textstyle{1\over 3}}\lambda^{1r}\tau^{st}\lambda_{rst}\}+O(n^{-1}).
\end{equation}

\section{Expectations of adjustments}

Detailed analysis given in the Appendix shows that we may approximate $E\{{\rm INF}(\psi)\}$ to $O(n^{-1})$ by \[
g_{{\rm INF}}(\theta)=\eta^{1/2}\lambda^{1r}\tau^{st}({\textstyle{1\over 2}}\lambda_{rs,t}+{\textstyle{1\over 6}}\lambda_{rst}),\]
 and
 $E\{{\rm NP}(\psi)\}$ to the same order by
 \[
 g_{{\rm NP}}(\theta)=-\eta^{1/2}\lambda^{1r}\nu^{st}(\lambda_{rs,t}+{\textstyle{1\over 2}}\lambda_{rst}).
  \]

These expansions permit a full statistical interpretation of the adjustment terms ${\rm NP}(\psi)$ and ${\rm INF}(\psi)$, which we do through a series of remarks.

\noindent{{\it Remark 1}. We begin by examining $E\{R(\psi)\}$ when there are no nuisance parameters. If nuisance parameters are absent, then $\lambda^{11}=(\lambda_{11})^{-1}$, $\eta=-\lambda_{11}$, $\tau^{11}=(-\lambda_{11})^{-1}$, and $\nu^{11}=0$, and it follows that
\[
E\{R(\psi)\}=(-\lambda_{11})^{-3/2}({\textstyle{1\over 2}}\lambda_{11,1}+{\textstyle{1\over 6}}\lambda_{111})+O(n^{-1}).
\]

\noindent{\it Remark 2}. The quantities $g_{\rm INF}(\theta)$ and $g_{\rm NP}(\theta)$ are related to asymptotic quantities detailed by Efron (1987) in description of the `bias corrected accelerated', $BC_a$, method of construction of bootstrap confidence intervals, which is analysed in detail by DiCiccio \& Efron (1996).
Specifically, we have $g_{\rm INF}(\theta)=a_0$ and $g_{\rm NP}(\theta)=z_0-a_0$, where $a_0=a_0(\theta)$ and $z_0=z_0(\theta)$ are respectively acceleration and bias-correction quantities. The quantity $a_0$ satisfies (DiCiccio \& Efron, 1996)
\[
a_0=-\frac{1}{6}\{{\rm skew}(U)+ {\rm skew}(T)\}+O(n^{-1}),
\]
where $U=(\hat \psi-\psi)/\sigma$, with $\sigma^2$ the variance of $\hat \psi$, given by
$\sigma^2 \equiv \sigma^2(\theta)=\lambda^{1,1}+ O(n^{-2})$, and $T=(\hat \psi-\psi)/\hat \sigma$, with $\hat \sigma^2=\sigma^2(\hat \theta)$. Further, $z_0$ is interpreted by
\[
\Phi(z_0)={\rm Pr}(\hat \psi \leq \psi)+O(n^{-1}),
\]
where $\Phi$ is the standard normal distribution function.

DiCiccio \& Efron (1996) note that the quantities $a_0$ and $z_0$ are invariant under reparameterisations of the model. Therefore, in using the asymptotic adjustment expectations $g_{\rm INF}(\theta)$ and $g_{\rm NP}(\theta)$ to interpret nuisance parameter effects on the inference on $\psi$, there is no restriction in assuming that the model under analysis is parameterised so that the interest parameter $\psi$ and the nuisance parameter $\phi$ are orthogonal (Cox \& Reid, 1987).
Therefore, now suppose there is a vector nuisance parameter $\phi$ present, but assume that the interest parameter $\psi$ and the nuisance parameter $\phi$ are orthogonal; then $\lambda^{11}=(\lambda_{11})^{-1}$, $\eta=-\lambda_{11}$, $\lambda^{1a}=0$ $(a=2,\ldots,d)$, $\tau^{rs}=0$ except when $r=s=1$, in which case $\tau^{11}=(-\lambda_{11})^{-1}$, and
\[
E\{{\rm INF}(\psi)\}=-(-\lambda_{11})^{-3/2}({\textstyle{1\over 2}}\lambda_{11,1}+{\textstyle{1\over 6}}\lambda_{111})+O(n^{-1}).
\]
Therefore, following Remark 1, to error of order $O(n^{-1})$, $E\{{\rm INF}(\psi)\}$ is seen to correspond to a mean adjustment for the signed root statistic $R(\psi)$ in the problem where the orthogonal nuisance parameter $\phi$ is known. Since the standard normal approximation to the distribution of $R(\psi)$ is typically rather accurate in scalar parameter cases without nuisance parameters, the mean adjustment should be quantitatively small quite generally, so we can anticipate that ${\rm INF}(\psi)$ is typically small.

\noindent{\it Remark 3}.
For general parameterisations, we have $\nu^{11}=\nu^{a1}=\nu^{1b}=0$ for $a,b=2,\ldots,d$, and thus,
\begin{eqnarray*}
E\{{\rm NP}(\psi)\}&=&-\eta^{1/2}\lambda^{1r}\nu^{ab}(\lambda_{ra,b}+{\textstyle{1\over 2}}\lambda_{rab})+O(n^{-1}) \\
&=&\eta^{1/2}\lambda^{1r}\nu^{ab}({\textstyle{1 \over 2}}\lambda_{rab}-\lambda_{ra/b})+O(n^{-1}),
\end{eqnarray*}
where $\lambda_{ra/b}=\partial\lambda_{ra}/\partial\theta^b$ and $\lambda_{ra/b}=\lambda_{ra,b}+\lambda_{rab}$.

Under orthogonality, $\nu^{ab}=\lambda^{ab}$ for $a,b=2,\ldots,d$, and the condition $\lambda_{1a}=0$ for $a=2,\ldots,d$ implies that $\lambda_{1a/b}=0$ for $b=2,\ldots,d$, so that the identity $\lambda_{1a/b}=\lambda_{1a,b}+\lambda_{1ab}$ yields $\lambda_{1a,b}=-\lambda_{1ab}$ for $a,b=2,\ldots,d$. Hence, nuisance parameter effects may be quantified from the expression $$E\{{\rm NP}(\psi)\}=-{\textstyle{1\over 2}}(-\lambda_{11})^{-1/2}\lambda^{ab}\lambda_{ab1}+O(n^{-1}).$$
Note that this gives $\beta_1=\eta^{1/2}E\{{\rm NP}(\psi)\}+O(n^{-1/2})=-\frac{1}{2}\lambda^{ab}\lambda_{ab1}+O(n^{-1/2})$.
Since the expansion for $E\{{\rm NP}(\psi)\}$ involves a multiple sum over the nuisance parameters, we see that ${\rm NP}(\psi)$ can be anticipated to be large when the number of nuisance parameters is large.

\noindent{\it Remark 4}.
Some further insight into ${\rm NP}(\psi)$ in the orthogonal case can be gleaned by noting that
\[
\frac{\partial\log\det[-L_{ab}\{\tilde\theta(\psi)\}]}{\partial\psi}=
L^{ab}(\theta)L_{ab1}(\theta)+O_p(n^{-1/2})=\lambda^{ab}\lambda_{ab1}+O_p(n^{-1/2}),
\]
which further relates $E\{{\rm NP}(\psi)\}$ to the specific adjustment function of Cox \& Reid (1987).
Thus, in this orthogonal case, if $\log\det\{-L_{ab}(\theta)\}$ does not change rapidly with $\psi$, such as when $L(\theta)=g(\psi)+h(\phi)$, in which case $\det\{-L_{ab}(\theta)\}$ is constant with respect to $\psi$, then $\lambda^{ab}\lambda_{ab1}$ is small in magnitude, and hence, we would expect ${\rm NP}(\psi)$ to be small in magnitude; see also the discussion in Cox \& Reid (1987).

\noindent{\it Remark 5}.
There is one further interpretation of ${\rm NP}(\psi)$ that is worth noting. DiCiccio \& Stern (1994a) showed that the difference between $\bar\psi$ and $\hat\psi$ is
\[
\bar\psi-\hat\psi=-\lambda^{11}\beta_1+O_p(n^{-3/2})=
\eta^{-1}\beta_1+O_p(n^{-3/2})=\eta^{-1/2}E\{NP(\psi)\}+O_p(n^{-3/2}),
\]
and hence, this difference, when in expressed in terms of standard deviations of $\hat\psi$, is
\[
\frac{\bar\psi-\hat\psi}{\eta^{-1/2}}=E\{NP(\psi)\}+O_p(n^{-1}).
\]

\noindent{\it Remark 6}.
 Note that the quantities $g_{{\rm NP}}(\theta)$ and $g_{{\rm INF}}(\theta)$ are both of order $O(n^{-1/2})$. As we shall illustrate, calculation of the individual values provides important statistical insight. We propose further that a simple measure of the relative influence within the assumed model of the nuisance parameter on inference on the interest parameter $\psi$, independent of the sample size $n$, might be obtained by considering their ratio $g_{{\rm NP}}(\theta)/g_{{\rm INF}}(\theta)$.

\noindent{\it Remark 7}. In general, the quantities $g_{{\rm NP}}(\theta)$ and $g_{{\rm INF}}(\theta)$ depend on the unknown parameter $\theta$.
In practice, following the bootstrap principle, they may be estimated by $g_{{\rm NP}}(\tilde \theta)$ and $g_{{\rm INF}}(\tilde \theta)$ respectively. An adjusted version of the signed root statistic $R(\psi)$, easily calculated in practice, once $g_{{\rm NP}}(\theta)$ and $g_{{\rm INF}}(\theta)$ have been calculated, is given by $R_a(\psi)=R(\psi)+g_{{\rm NP}}(\tilde \theta)+g_{{\rm INF}}(\tilde \theta)$. Since $g_{{\rm NP}}(\tilde \theta)-g_{{\rm NP}}(\theta)=O_p(n^{-1})$, we have that $R_a(\psi)=R^*(\psi)+O_p(n^{-1})$, and therefore that $R_a(\psi)$ has the standard normal distribution to error of order $O(n^{-1})$. DiCiccio \& Efron (1996) previously remarked that $R(\psi)+z_0(\hat \theta)$ is standard normal to error of order $O(n^{-1})$, but did not investigate practical use of this statistic for inference: an alternative is the statistic $R_a(\psi)=R(\psi)+z_0(\tilde \theta)$. Although no claim of desirable large deviation properties of the kind enjoyed by the method of Skovgaard (1986) can be made for this statistic, empirical evidence, not reported here, suggests that it nevertheless yields highly accurate inference in many settings.

\noindent{\it Remark 8}.
Note that the asymptotic regime adopted here is one in which the dimensionality $d-1$ of the nuisance parameter $\phi$ remains fixed as the sample size $n$ increases. However, we propose that examination of the quantities
$g_{{\rm NP}}(\theta)$ and $g_{{\rm INF}}(\theta)$ and their ratio
is a useful device to quantify the effect of an increasing dimension of nuisance parameter on the inference, as we shall illustrate in the next Section. For stratified models, such as those in Examples 2, 4, 5 and 6 below, Sartori
(2003) noted that, when both the sample size $n$ within each stratum and
the number of nuisance parameters $q$ tend to infinity, $NP( \psi) = O_p(q m^{-1/2})$, while
$INF( \psi) = O_p(m^{-1/2})$, where $m=nq$ is the total sample size, irrespective of the nature of the sequence $\{q,n\}$. Hence, the ratio $NP(\psi)/INF( \psi) = O_p(q)$ in such an asymptotic regime,
consistent with calculations given in
Examples 2, 4, 5 and 6 below. Relative to the inference adjustment, the nuisance parameter adjustment increases at a rate proportional to the dimension of the nuisance parameter.

\section{Examples}

We consider here a number of theoretical and numerical examples.

Example 1. {\it Normal linear regression.} Let $Y_1, \ldots, Y_n$ denote independent random variables of the form $Y_i=x_i^T\beta+\sigma \epsilon_i$, where $x_1, \ldots, x_n$ are known covariate vectors of length $q$, $\sigma$ is an unknown scalar interest parameter and $\beta$ is an unknown nuisance parameter vector of length $q$, so that $\theta=(\sigma, \beta)$. The $\epsilon_i$ are assumed to be independent standard normal random variables.

In this case, $n^{1/2}g_{\rm INF}(\theta)=2^{1/2}/3$ and $n^{1/2}g_{\rm NP}(\theta)=q/2^{1/2}$. Note that these quantities do not depend on the parameter value $\theta$, while $\eta=2n/\sigma^2$. Nuisance parameter effects are determined, to second order, only by the dimensionality of the nuisance parameter $\beta$, not its value. This observation in turn would suggest that inference based on the bootstrap distribution of $R(\sigma)$ should be highly accurate. In fact, $R(\sigma)$ is a simple function of $\hat \sigma^2/\sigma^2$, which has a distribution free of $\theta$: $(n-q) \hat \sigma^2/\sigma^2$ is distributed as chi-squared on $n-q$ degrees of freedom. A bootstrap calculation will, modulo simulation variability, reproduce the exact sampling distribution of $R(\sigma)$.

Example 2. {\it Neyman-Scott model.} Let $Y_{ij}$, for $i=1, \ldots , n$ and $j=1, \ldots , q$ be independent Gaussian random variables, with $Y_{ij}$ being distributed as $N(\mu_j, \sigma^2)$. The interest parameter is $\sigma$, with nuisance parameter $(\mu_1, \ldots, \mu_q)$, so that $\theta=(\sigma, \mu_1, \ldots, \mu_q)$.

Now we calculate $n^{1/2}g_{{\rm INF}}(\theta)=1/\{$1.5$ (2q)^{1/2}\}$, with $n^{1/2}g_{{\rm NP}}(\theta)= (q/2)^{1/2}$, so that $g_{{\rm NP}}(\theta)/g_{{\rm INF}}(\theta) = $1$\cdot$5$q$. Again, these quantities do not depend on the value of $\theta$, only the dimension $q$ of the nuisance parameter. The adjusted information is given by $\eta=2nq/\sigma^2$. As in Example 1, the signed root statistic $R(\sigma)$ has a distribution free of the parameter value: it is a function of the pivotal quantity $\hat \sigma^2/\sigma^2$, and its exact sampling distribution can be constructed by bootstrapping.

A related problem concerns a generalisation of the Behrens-Fisher problem, in which we observe $Y_{ij}$, for $i=1, \ldots , n$ and $j=1, \ldots , q$ to be independent Gaussian random variables, with $Y_{ij}$ being distributed as $N(\mu, \sigma_j^2)$. The interest parameter is the common mean $\mu$, with $(\sigma^2_1, \ldots, \sigma_q^2)$ as nuisance. In this case, we see that $E\{{\rm INF}(\psi)\}$ and $E\{{\rm NP}(\psi)\}$ are both $O(n^{-1})$, not $O(n^{-1/2})$. Nuisance parameter effects are quantitatively slight though, by contrast with what is noted above, in this case the signed root statistic $R(\mu)$ is not exactly pivotal, and the bootstrap inference is not exact. Limited numerical results given by Young (2009) for the case $q=2$ would indicate, however, that the bootstrap inference is highly accurate even for small sample size $n$.

Example 3. {\it Exponential regression.}  Suppose $Y_1, \ldots, Y_n$ are independent exponential random variables, with means depending on given covariate values. We suppose for simplicity the case of two covariates, though our conclusions extend immediately to the case with a general number of covariates. So, we suppose $Y_i$ is exponentially distributed with mean
$\phi_1 \exp(-\psi z_i-\phi_2 w_i)$, with $\sum z_i =\sum w_i=0$, and $\psi$ the interest parameter. Routine calculations show that
$g_{\rm INF}(\theta)$ and $g_{\rm NP}(\theta)$, though complicated functions of the covariate values, are again free of the parameter $\theta=(\psi, \phi_1, \phi_2)$. Further, the signed root statistic $R(\psi)$ is again easily seen to be exactly pivotal, and bootstrap inference is once more exact.

In the simple case of a single covariate, with $E(Y_i)=\phi \exp(-\psi z_i)$, with $\sum z_i=0$, we have
\[
E\{{\rm NP}(\psi)\}=0+O(n^{-1}), \;E\{{\rm INF}(\psi)\}=-(\sum z_i^2)^{-3/2} ({\textstyle{\frac{1}{6} }}\sum z_i^3)+O(n^{-1}):
\]
the nuisance parameter adjustment has expectation of smaller order of magnitude than that of the information adjustment.

We consider now from a numerical perspective three examples with many nuisance parameters previously discussed by Sartori et al.\ (1999). In each, we provide illustration of dependence of the measure $g_{{\rm NP}}(\theta)/g_{{\rm INF}}(\theta)$ on the dimensionality of the nuisance parameter.

Example 4. {\it Inverse Gaussian model.}  Let $Y_{ij}$, for $i=1, \ldots , n$ and $j=1, \ldots , q$ be independent, inverse Gaussian random variables, with $Y_{ij}$ having probability density
\[
f(y;\psi,\phi_j)=\{\psi/(2\pi)\}^{1/2} y^{-3/2} \exp \{-\textstyle\frac{1}{2}(\psi y^{-1}+\phi_j y)+(\psi
\phi_j)^{1/2}\}, \; \; y>0,
\]
where $\psi>0$ and $\phi_j>0$, so that $\theta=(\psi, \phi_1, \ldots, \phi_q)$ and the overall sample size is $m=nq$.

Simple algebraic manipulations show that, independently of the parameter value $\theta$, $n^{1/2}g_{{\rm INF}}(\theta) = -1/\{$1.5$ (2q)^{1/2}\}$, and
$n^{1/2}g_{{\rm NP}}(\theta)= -(q/2)^{1/2}$, so that \break $g_{{\rm NP}}(\theta)/g_{{\rm INF}}(\theta) = $1.5$q$ in this model. We note that in this model the adjusted information for $\psi$ is given by $\eta=nq/(2\psi^2)$.

Example 5. {\it Multi-sample exponential model.}  Let $Y_{ij}$, for $i=1, \ldots , n$ and $j=1, \ldots , q$ be independent, exponential random variables, with $Y_{ij}$ having mean $1/\phi_j$. The parameter of interest is
\[
\psi=q^{-1}\sum_{j=1}^q \exp (-\phi_j t_0),
\]
where $t_0 >0$ is a fixed constant and $\theta=(\psi, \phi)$, with the nuisance parameter $\phi=(\phi_2, \ldots, \phi_q)$. As noted by Sartori et al.\  (1999), $q\psi$ may be interpreted as the expected number of items failing by $t_0$ in a parallel system with failures rates $\phi_1, \ldots, \phi_q$.

The interest parameter $\psi$ is therefore a nonlinear function of the canonical parameter in a full exponential family model. Again, construction of the information and nuisance parameter adjustments ${\rm INF}(\psi)$ and ${\rm NP}(\psi)$ is straightforward, though the constrained maximum likelihood estimator $\tilde \theta$ must be calculated numerically.

By contrast with previous examples, in this model the ratio $g_{{\rm NP}}(\theta)/g_{{\rm INF}}(\theta)$ depends on the value of the parameter $\theta$. Values illustrating the effect of increasing nuisance parameter dimension are given in Table 1 for two cases. In both $t_0$=0.5: case (a) considers $\phi_i=1, i=1, \ldots, q$, so that $\psi=$0.6065; case (b) fixes $\psi=$0.0333 for each dimension of nuisance parameter, sets $\exp(-\phi_qt_0)=q\psi/2$ and fixes $\phi_1= \ldots = \phi_{q-1}$, the common value being determined by the specified $\psi$. Acute dependence of the ratio on the actual parameter values, rather than just the nuisance parameter dimension as in previous examples, is apparent.

\begin{table}
\caption{Dependence of ratio $g_{{\rm NP}}(\theta)/g_{{\rm INF}}(\theta)$ on $q$,
multi-sample exponential model. Case (a) has $\phi_i=1, i=1, \ldots , q$, case (b) has $\phi_1= \ldots = \phi_{q-1}$, with $\exp(-\phi_qt_0)=q\psi/2$. \vspace{0.25cm}}
\centering
\begin{tabular}{lrrrrr}
\vspace{0.2cm}$q$ & 2& 5& 10& 20& 50 \\
(a)& 2.25&9.00&20.25&42.75& 110.25\\
(b) &-2.10 &-5.50&-8.56 &-15.76 &-130.29  \\
\end{tabular}
\end{table}

Example 6. {\it Curved exponential family model.} Our final example concerns a model for which calculation of $R^*(\psi)$ is intractable: the sample space derivatives, derivatives of the log-likelihood with respect to the maximum likelihood estimator, required by the construction (\ref{definitioncorrection}) of $R^*(\psi)$, must be approximated. By contrast, the calculations required to evaluate $g_{\rm INF}(\theta)$ and $g_{\rm NP}(\theta)$ are no more complex than in the other examples.

Let $Y_{ij}$, for $i=1, \ldots , n$ and $j=1, \ldots , q$ be independent normal random variables with means $\mu_j>0$ and variances $\psi \mu_j^{1/2}$. This model constitutes a curved exponential family. The parameter of interest is $\psi$, with $\mu_1, \ldots, \mu_q$ as nuisance parameters, $\theta=(\psi, \mu_1, \ldots, \mu_q)$.

\begin{table}
\caption{Dependence of ratio $g_{{\rm NP}}(\theta)/g_{{\rm INF}}(\theta)$ on $q$,
multi-sample curved exponential family model. Case (a) has $\psi=1, \mu_i=i, i=1, \ldots , q$, case (b) has $\psi=1, \mu_i=1, i=1, \ldots, q$. \vspace{0.25cm}}
\centering
\begin{tabular}{lrrrrrr}
\vspace{0.2cm}$q$ &1& 2& 5& 10& 20& 50 \\
(a)& 1.11&2.45&6.77&14.17&29.09&74.01\\
(b) &1.11&2.21&5.53&11.05&22.11&55.26  \\
\end{tabular}

\end{table}

Again, the ratio $g_{{\rm NP}}(\theta)/g_{{\rm INF}}(\theta)$ depends on the value of the parameter $\theta$. Illustrative values are given in Table 2, for two cases: case (a) has $\psi=1, \mu_i=i, i=1, \ldots, q$, while case (b) has $\psi=1, \mu_i=1, i=1, \ldots, q$.

\section{Decomposition of the Bartlett correction factor}

Recall that the sum of $g_{\rm INF}(\theta)$ and $g_{\rm NP}(\theta)$ is, to $O(n^{-1})$, equal to 
\begin{align*}
E\{-R(\psi)\}&=-\eta^{1/2}(\lambda^{1r}\lambda^{st}\lambda_{rs,t}
+{\textstyle{1\over 2}}\lambda^{1r}\tau^{st}\lambda_{rs,t}
+{\textstyle{1\over 2}}\lambda^{1r}\lambda^{st}\lambda_{rst}
+{\textstyle{1\over 3}}\lambda^{1r}\tau^{st}\lambda_{rst})\\
&=-\eta^{1/2}\lambda^{1r}\lambda^{st}(\lambda_{rs,t}+{\textstyle{1\over 2}}\lambda_{rst})
-\eta^{1/2}\lambda^{1r}\tau^{st}({\textstyle{1\over 2}}\lambda_{rs,t}+{\textstyle{1\over 3}}\lambda_{rst}).
\end{align*}

To decide how we might choose $g_{\rm INF}(\theta)$ and $g_{\rm NP}(\theta)$ in a decomposition of this sum, consider imposing two conditions: first, $g_{\rm INF}(\theta)$ must take the same value whether we have no nuisance parameters or we have orthogonal nuisance parameters; and second, $g_{\rm NP}(\theta)$ must be $0$ when we have no nuisance parameters. These conditions suggest that $\tau^{rs}$ and $\nu^{rs}$ play a key role. Note that $\tau^{11}=(-\lambda_{11})^{-1}$ when there are no nuisance parameters, while for orthogonal nuisance parameters $\tau^{rs}=0$ except when $r=s=1$, in which case $\tau^{11}=(-\lambda_{11})^{-1}$. Thus, $\tau^{rs}$ is the same in the orthogonal nuisance parameter case as it is when nuisance parameters are absent. On the other hand, since $\nu^{rs}=0$ whenever either or both of $r$ and $s$ are $1$, we have that $\nu^{11}=0$ when there are no nuisance parameters. It is readily seen that the decomposition of the sum into $g_{\rm INF}(\theta)$ and $g_{\rm NP}(\theta)$ according to the two conditions can be achieved if we substitute $\lambda^{st}=\nu^{st}-\tau^{st}$ in the sum and then take $g_{\rm INF}(\theta)$ to consist of those terms involving $\tau^{st}$ and take $g_{\rm NP}(\theta)$ to consist of those terms involving $\nu^{st}$. We demonstrate here that the same reasoning may be applied to obtain a decomposition of the Bartlett correction factor for the likelihood ratio statistic $W(\psi)$.

Lawley (1956) showed (see also DiCiccio \& Stern, 1994a) that the expectation of $W(\psi)$ is $E\{W(\psi)\}=1+b(\theta)+O(n^{-3/2})$, where 
\begin{align*}
b(\theta)&=(\lambda^{rs}\lambda^{tu}-\nu^{rs}\nu^{tu})({\textstyle{1\over 4}}\lambda_{rstu}-\lambda_{rst/u}+\lambda_{rt/su})\\
&\qquad-(\lambda^{rs}\lambda^{tu}\lambda^{vw}-\nu^{rs}\nu^{tu}\nu^{vw})({\textstyle{1\over 4}}\lambda_{rst}\lambda_{uvw}-\lambda_{rst}\lambda_{uv/w}+\lambda_{rs/t}\lambda_{uv/w})\\
&\qquad\qquad-(\lambda^{ru}\lambda^{sw}\lambda^{tv}-\nu^{ru}\nu^{sw}\nu^{tv})({\textstyle{1\over 6}}\lambda_{rst}\lambda_{uvw}-\lambda_{rst}\lambda_{uv/w}+\lambda_{rs/t}\lambda_{uv/w}).
\end{align*}

We now decompose $b(\theta)$ into the sum $b(\theta)=b_{\rm INF}(\theta)+b_{\rm NP}(\theta)$, where $b_{\rm INF}(\theta)$ is the same whether we have no nuisance parameters or whether we have orthogonal nuisance parameters, and $b_{\rm NP}(\theta)$ is $0$ when there are no nuisance parameters. We make the substitution $\lambda^{rs}=\nu^{rs}-\tau^{rs}$ in $b(\theta)$: $b_{\rm INF}(\theta)$ consists of those terms involving the $\tau^{rs}$ but not the $\nu^{rs}$; $b_{\rm NP}(\theta)$ consists of those terms that involve the $\nu^{rs}$ in any way. 

Succinct expressions for $b_{\rm INF}(\theta)$ and $b_{\rm NP}(\theta)$ derived this way are
\begin{align*}
b_{\rm INF}(\theta)&=\tau^{rs}\tau^{tu}({\textstyle{1\over 4}}\lambda_{rstu}-\lambda_{rst/u}+\lambda_{rt/su})\\
&\qquad+\tau^{rs}\tau^{tu}\tau^{vw}({\textstyle{1\over 4}}\lambda_{rst}\lambda_{uvw}-\lambda_{rst}\lambda_{uv/w}+\lambda_{rs/t}\lambda_{uv/w})\\
&\qquad\qquad+\tau^{ru}\tau^{sw}\tau^{tv}({\textstyle{1\over 6}}\lambda_{rst}\lambda_{uvw}-\lambda_{rst}\lambda_{uv/w}+\lambda_{rs/t}\lambda_{uv/w}),
\end{align*}
and
\begin{align*}
b_{\rm NP}(\theta)&=(\lambda^{rs}\lambda^{tu}-\tau^{rs}\tau^{tu}-\nu^{rs}\nu^{tu})({\textstyle{1\over 4}}\lambda_{rstu}-\lambda_{rst/u}+\lambda_{rt/su})\\
&-(\lambda^{rs}\lambda^{tu}\lambda^{vw}+\tau^{rs}\tau^{tu}\tau^{vw}-\nu^{rs}\nu^{tu}\nu^{vw})({\textstyle{1\over 4}}\lambda_{rst}\lambda_{uvw}-\lambda_{rst}\lambda_{uv/w}+\lambda_{rs/t}\lambda_{uv/w})\\
&-(\lambda^{ru}\lambda^{sw}\lambda^{tv}+\tau^{ru}\tau^{sw}\tau^{tv}-\nu^{ru}\nu^{sw}\nu^{tv})({\textstyle{1\over 6}}\lambda_{rst}\lambda_{uvw}-\lambda_{rst}\lambda_{uv/w}+\lambda_{rs/t}\lambda_{uv/w}).
\end{align*}

If there are no nuisance parameters or there are orthogonal nuisance parameters, then 
\begin{align*}
b_{\rm INF}(\theta)&=(\lambda_{11})^{-2}({\textstyle{1\over 4}}\lambda_{1111}-\lambda_{111/1}+\lambda_{11/11})\\
&\qquad-(\lambda_{11})^{-3}({\textstyle{1\over 4}}\lambda_{111}\lambda_{111}-\lambda_{111}\lambda_{11/1}+\lambda_{11/1}\lambda_{11/1})\\
&\qquad\qquad-(\lambda_{11})^{-3}({\textstyle{1\over 6}}\lambda_{111}\lambda_{111}-\lambda_{111}\lambda_{11/1}+\lambda_{11/1}\lambda_{11/1}).
\end{align*}

Note that if there are no nuisance parameters, $\nu^{11}=0$ and $\tau^{11}=-\lambda^{11}$, so that $b_{\rm NP}(\theta)$ is identically zero. It is useful to evaluate $b_{\rm NP}(\theta)$ in the case of orthogonal nuisance parameters to show better the effect of nuisance parameters. Now, by making the substitution $\lambda^{rs}=\nu^{rs}-\tau^{rs}$, we have 
\begin{align*}
b_{\rm NP}(\theta)&=\{(\nu^{rs}-\tau^{rs})(\nu^{tu}-\tau^{tu})-\tau^{rs}\tau^{tu}-\nu^{rs}\nu^{tu}\}({\textstyle{1\over 4}}\lambda_{rstu}-\lambda_{rst/u}+\lambda_{rt/su})\\
&\qquad-\{(\nu^{rs}-\tau^{rs})(\nu^{tu}-\tau^{tu})(\nu^{vw}-\tau^{vw})+\tau^{rs}\tau^{tu}\tau^{vw}-\nu^{rs}\nu^{tu}\nu^{vw}\}\\
&\qquad\quad\quad\quad\quad\times({\textstyle{1\over 4}}\lambda_{rst}\lambda_{uvw}-\lambda_{rst}\lambda_{uv/w}+\lambda_{rs/t}\lambda_{uv/w})\\
&\qquad\qquad-\{(\nu^{ru}-\tau^{ru})(\nu^{sw}-\tau^{sw})(\nu^{tv}-\tau^{tv})+\tau^{ru}\tau^{sw}\tau^{tv}-\nu^{ru}\nu^{sw}\nu^{tv}\}\\
&\qquad\qquad\quad\quad\quad\times({\textstyle{1\over 6}}\lambda_{rst}\lambda_{uvw}-\lambda_{rst}\lambda_{uv/w}+\lambda_{rs/t}\lambda_{uv/w})\\
\end{align*}
\begin{align*}
&=-(\tau^{rs}\nu^{tu}+\nu^{rs}\tau^{tu})({\textstyle{1\over 4}}\lambda_{rstu}-\lambda_{rst/u}+\lambda_{rt/su})\\
&\qquad-(\tau^{rs}\tau^{tu}\nu^{vw}+\tau^{rs}\nu^{tu}\tau^{vw}+\nu^{rs}\tau^{tu}\tau^{vw}-\tau^{rs}\nu^{tu}\nu^{vw}-\nu^{rs}\tau^{tu}\nu^{vw}-\nu^{rs}\nu^{tu}\tau^{vw})\\
&\qquad\quad\quad\quad\quad\times({\textstyle{1\over 4}}\lambda_{rst}\lambda_{uvw}-\lambda_{rst}\lambda_{uv/w}+\lambda_{rs/t}\lambda_{uv/w})\\
&\qquad\qquad-(\tau^{ru}\tau^{sw}\nu^{tv}+\tau^{ru}\nu^{sw}\tau^{tv}+\nu^{ru}\tau^{sw}\tau^{tv}-\tau^{ru}\nu^{sw}\nu^{tv}-\nu^{ru}\tau^{sw}\nu^{tv}-\nu^{ru}\nu^{sw}\tau^{tv})\\
&\qquad\qquad\quad\quad\quad\quad\times({\textstyle{1\over 6}}\lambda_{rst}\lambda_{uvw}-\lambda_{rst}\lambda_{uv/w}+\lambda_{rs/t}\lambda_{uv/w}).
\end{align*}

We consider each of the terms in $b_{\rm NP}(\theta)$ separately under orthogonality:
\begin{align*}
&-(\tau^{rs}\nu^{tu}+\nu^{rs}\tau^{tu})({\textstyle{1\over 4}}\lambda_{rstu}-\lambda_{rst/u}+\lambda_{rt/su})\\
&\qquad =
(\lambda_{11})^{-1}\lambda^{ab}({\textstyle{1\over 2}}\lambda_{11ab}-\lambda_{1ab/1}-\lambda_{11a/b});\\ \\ 
&-(\tau^{rs}\tau^{tu}\nu^{vw}+\tau^{rs}\nu^{tu}\tau^{vw}+\nu^{rs}\tau^{tu}\tau^{vw}-\tau^{rs}\nu^{tu}\nu^{vw}-\nu^{rs}\tau^{tu}\nu^{vw}-\nu^{rs}\nu^{tu}\tau^{vw})\\ 
&\quad\quad\quad\times({\textstyle{1\over 4}}\lambda_{rst}\lambda_{uvw}-\lambda_{rst}\lambda_{uv/w}+\lambda_{rs/t}\lambda_{uv/w})\\
&\qquad=-(\lambda_{11})^{-2}\lambda^{ab}({\textstyle{1\over 2}}\lambda_{111}\lambda_{1ab}+{\textstyle{1\over 4}}\lambda_{11a}\lambda_{11b}-\lambda_{1ab}\lambda_{11/1}+
\lambda_{11/1}\lambda_{ab/1})\\
&\qquad\quad\quad\quad-(\lambda_{11})^{-1}\lambda^{ab}\lambda^{cd}({\textstyle{1\over 2}}\lambda_{11a}\lambda_{bcd}+{\textstyle{1\over 4}}\lambda_{1ab}\lambda_{1cd}-\lambda_{11a}\lambda_{bc/d}-
\lambda_{1ab}\lambda_{1c/d}+\lambda_{11/a}\lambda_{bc/d});\\ \\
&-(\tau^{ru}\tau^{sw}\nu^{tv}+\tau^{ru}\nu^{sw}\tau^{tv}+
\nu^{ru}\tau^{sw}\tau^{tv}-\tau^{ru}\nu^{sw}\nu^{tv}-\nu^{ru}\tau^{sw}\nu^{tv}-
\nu^{ru}\nu^{sw}\tau^{tv})\\
&\quad\quad\quad\times({\textstyle{1\over 6}}\lambda_{rst}\lambda_{uvw}-\lambda_{rst}\lambda_{uv/w}+\lambda_{rs/t}\lambda_{uv/w})\\
&\qquad=-(\lambda_{11})^{-2}\lambda^{ab}({\textstyle{1\over 2}}\lambda_{11a}\lambda_{11b}-\lambda_{11a}\lambda_{11/b}-\lambda_{11a}\lambda_{1b/1})\\
&\qquad\quad\quad\quad-(\lambda_{11})^{-1}\lambda^{ab}\lambda^{cd}({\textstyle{1\over 2}}\lambda_{1ac}\lambda_{1bd}-\lambda_{1ac}\lambda_{bd/1}).
\end{align*}

The resulting formula for $b_{\rm NP}(\theta)$ in the presence of orthogonal nuisance parameters is
\begin{align*}
b_{\rm NP}(\theta)&=(\lambda_{11})^{-1}\lambda^{ab}({\textstyle{1\over 2}}\lambda_{11ab}-\lambda_{1ab/1}-\lambda_{11a/b})\\
&\qquad-(\lambda_{11})^{-2}\lambda^{ab}(\lambda_{111}\lambda_{1ab}+{\textstyle{3\over 4}}\lambda_{11a}\lambda_{11b}\\
&\qquad\quad\quad\quad\quad\quad-\lambda_{11a}\lambda_{11/b}-\lambda_{11a}\lambda_{1b/1}-
\lambda_{1ab}\lambda_{11/1}+\lambda_{11/1}\lambda_{ab/1})\\
&\qquad\qquad-(\lambda_{11})^{-1}\lambda^{ab}\lambda^{cd}({\textstyle{1\over 2}}\lambda_{11a}\lambda_{bcd}+{\textstyle{1\over 4}}\lambda_{1ab}\lambda_{1cd}+{\textstyle{1\over 2}}\lambda_{1ac}\lambda_{1bd}\\
&\qquad\qquad\quad\quad\quad\quad\quad-\lambda_{11a}\lambda_{bc/d}-\lambda_{1ab}\lambda_{1c/d}-
\lambda_{1ac}\lambda_{bd/1}+\lambda_{11/a}\lambda_{bc/d}.
\end{align*}
Just as for $g_{\rm NP}(\theta)$ in the case of orthogonal nuisance parameters, we see that $b_{\rm NP}(\theta)$ involves multiple sums over the indices for the nuisance parameters, so $b_{\rm NP}(\theta)$ can be expected to be large when the number of nuisance parameters is large. 

An interesting feature emerges from comparing the formulas for $g_{\rm NP}(\theta)$ and $b_{\rm NP}(\theta)$ in the orthogonal nuisance parameter case. While the expression for $g_{\rm NP}(\theta)$ involves a double sum over the indices for the nuisance parameters, the expression for $b_{\rm NP}(\theta)$ involves both double and quadruple sums. Consequently, we might reasonably expect the ratio $b_{\rm NP}(\theta)/b_{\rm INF}(\theta)$ to grow more rapidly with the number of nuisance parameters than does the ratio $g_{\rm NP}(\theta)/g_{\rm INF}(\theta)$. This phenomenon is apparent in Example 1, for which $g_{\rm NP}(\theta)/g_{\rm INF}(\theta)=3q/2.$ It turns out that $b_{\rm INF}(\theta)=n^{-1}\textstyle{1 \over 3}$ and $b_{\rm NP}(\theta)=n^{-1}(q^2+q)$, so $b_{\rm NP}(\theta)/b_{\rm INF}(\theta)=3(q^2+q)$. In this example, the ratio $b_{\rm NP}(\theta)/b_{\rm INF}(\theta)$ grows quadratically with the number of nuisance parameters, while the ratio $g_{\rm NP}(\theta)/g_{\rm INF}(\theta)$ only grows linearly. 

\section{Decompositions for other pivots}

So far, our focus has been on inference based on an adjusted version of the signed root likelihood ratio statistic; however, other pivots that are asymptotically standard normal also find widespread use, notably the Wald-type pivots based on the difference $\hat\psi-\psi$ and the score-type pivots based on the derivative $M_1(\psi)=d M(\psi)/d\psi=L_1\{\tilde\theta(\psi)\}$.
DiCiccio et al.\  (2015) provide analysis of circumstances where inference, such as $p-$values, obtained by bootstrapping various first-order asymptotically equivalent pivots will agree to higher-order with that obtained from the signed root statistic. It is of interest to assess the impact that nuisance parameters have on higher-order adjustments obtained by Cornish-Fisher transformation to these other pivots.  We examine the structure of these adjustments in terms of the quantities $g_{\rm INF}(\theta)$ and $g_{\rm NP}(\theta)$, to allow explicit comparisons with inference based on $R(\psi)$.

Let $T(\psi)$ denote an asymptotically standard normal pivot, and let its cumulants be denoted by $\kappa_1$, $\kappa_2$, etc. Typically, the mean $\kappa_1$ and  skewness $\kappa_3$ are of order $O(n^{-1/2})$, while the variance $\kappa_2=1+O(n^{-1})$; the fourth and higher-order cumulants are of order $O(n^{-1})$ or smaller. Central to higher-order inference based on $T(\psi)$ is the Cornish-Fisher transformation $T-{1\over 6}\kappa_3T^2-\kappa_1+{1 \over 6}\kappa_3$, which has the standard normal distribution to error of order $O(n^{-1})$. The Cornish-Fisher transformation of $R(\psi)$ agrees with the $R^*(\psi)$ statistic to error of order $O(n^{-1})$. The adjustment terms ${1\over 6}\kappa_3$ and $-\kappa_1+{1 \over 6}\kappa_3$ that appear in the Cornish-Fisher transformation depend on $\theta$, so they would need to be estimated to achieve higher-order inference in practice. An interpretation of the adjustment made by the Cornish-Fisher transformation is that whether or not a mean adjustment suffices to make the desired correction hinges on the order of $\kappa_{3}$. This is an important factor differentiating the signed root statistic from other asymptotically standard normal pivots.

We report $\kappa_1$ and $\kappa_3$ for some common choices of $T(\psi)$. For $T(\psi)=R(\psi)$, we have seen that $\kappa_1=-g_{\rm INF}(\theta)-g_{\rm NP}(\theta)+O(n^{-1})$; in this case, $\kappa_3=O(n^{-1})$. Consequently, higher-order inference based on $R(\psi)$ requires estimation of $\kappa_1$ only, and estimation of $\kappa_3$ is not necessary. 

To report $\kappa_1$ and $\kappa_3$ for other pivots $T(\psi)$, it is convenient to introduce one further asymptotic quantity in addition to $g_{\rm INF}(\theta)$ and $g_{\rm NP}(\theta)$. This quantity is $d \equiv d(\theta)=-\eta^{1/2}{1 \over 6}\lambda^{1r}\tau^{st}\lambda_{rst}$, which arises quite naturally from the profile log-likelihood function. It turns out that the third derivative of the profile log-likelihood function evaluated at $\hat\psi$ is $M_3(\hat\psi)=\eta^{3/2}6d+O_p(n^{1/2})$. The quantity $d$ is also related to Efron's (1987) asymptotic adjustments $a_0$ and $c_q$, which were discussed by DiCiccio \& Efron (1996): $d=2a_0+c_q$. Furthermore, in terms of $g_{\rm INF}(\theta)$, $g_{\rm NP}(\theta)$, and $d$, the mean of $\hat\psi$ is $E(\hat\psi)=\psi-(2g_{\rm INF}(\theta)+g_{\rm NP}(\theta)-d)\eta^{-1/2}+O(n^{-3/2})$.

A key property of the quantity $d$ is that it is the same whether there are no nuisance parameters or there are orthogonal nuisance parameters. In both cases, the formula for $d$ becomes $d=-(-\lambda_{11})^{-3/2}{1 \over 6}\lambda_{111}$. Thus, $d$ is similar to $g_{\rm INF}(\theta)$: we would not expect $d$ to grow with the number of nuisance parameters. The quantity $d$ does differ from $g_{\rm INF}(\theta)$ and $g_{\rm NP}(\theta)$ in one important respect: while $g_{\rm INF}(\theta)$ and $g_{\rm NP}(\theta)$ are invariant under reparameterizations $\theta=(\psi,\phi)\rightarrow \{g(\psi), h(\psi,\phi)\}$, where $\phi=(\theta^2,\ldots,\theta^d)$ contains the nuisance parameters and $g(\psi)$ is a monotonically increasing function, $d$ does not enjoy the property of invariance. 

We next consider the Wald statistic with observed information, $T(\psi)=(\hat\psi-\psi)/(-\hat L^{11})^{1/2}$, and the Wald statistic with expected information, $T(\psi)=(\hat\psi-\psi)/(-\hat\lambda^{11})^{1/2}=(\hat\psi-\psi)\hat\eta^{1/2}$. The distributions of these pivots are the same to error of order $O(n^{-1})$. For both Wald statistics, $\kappa_1=-\{g_{\rm INF}(\theta)+g_{\rm NP}(\theta)+d\}+O(n^{-1})$ and $\kappa_3=-6d+O(n^{-1})$. Consequently, the Wald statistics are similar to the signed root of the likelihood ratio statistic in that nuisance parameters affect the higher-order adjustment terms through $g_{\rm NP}(\theta)$, which is involved in $\kappa_1$.

Finally, we consider the score statistic with observed information, $T(\psi)=M_1(\psi)(-\hat L^{11})^{1/2}$, and the score statistic with expected information, $T(\psi)=M_1(\psi)(-\hat\lambda^{11})^{1/2}=M_1(\psi)\hat\eta^{-1/2}$. Just as for the Wald statistics discussed above, the distributions of these pivots agree to error of order $O(n^{-1})$; for these score statistics, $\kappa_1=-\{g_{\rm INF}(\theta)+g_{\rm NP}(\theta)-2d\}+O(n^{-1})$ and $\kappa_3=12d+O(n^{-1})$. Again, nuisance parameters influence the higher-order adjustment terms through $g_{\rm NP}(\theta)$, which is a component of $\kappa_1$. 

An important property of the profile log-likelihood function $M(\psi)$ is that the expectation of the profile score is $E\{M_{1}(\psi)\}=-\eta^{1/2}g_{\rm NP}(\theta)+O(n^{-1})$. Thus, $E\{M_{1}(\psi)\}$ is of order $O(1)$; the expectation of the profile score does even vanish asymptotically. Adjusted profile likelihood is discussed in the Appendix. Most of the adjustment functions $B(\psi)$ that have been proposed to construct an adjusted profile log-likelihood $\bar M(\psi)=M(\psi)+B(\psi)$ have the property that $E\{B_1(\psi)\}=\eta^{1/2}g_{\rm NP}(\theta)+O(n^{-1})$, so the expectation of the adjusted profile score is $E\{M_1(\psi)\}=O(n^{-1})$, which does vanish asymptotically. 

For $T(\psi)=\bar R(\psi)={\rm sgn}(\bar \psi-\psi)[2\{ \bar M(\bar\psi)- \bar M(\psi) \}]^{1/2}$, as detailed in the Appendix, we have $\kappa_1=-g_{\rm INF}(\theta)+O(n^{-1})$ and $\kappa_3=O(n^{-1})$. Thus, at order $O(n^{-1/2})$, the difference between the distribution of $\bar R(\psi)$ and the standard normal distribution depends on $g_{\rm INF}(\theta)$, a term which is the same whether there are no nuisance parameters present or there are orthogonal nuisance parameters. Consequently, we expect the difference between the distribution of $\bar R(\psi)$ and the standard normal distribution not to grow inordinately as the number of nuisance parameters increase. 

Similar comments apply to Wald statistics and score statistics based on the adjusted profile log-likelihood function. For example, for $T(\psi)=(\bar \psi-\psi)\{-\bar M_{11}(\bar\psi)\}^{1/2}$, we have $\kappa_1=-\{g_{\rm INF}(\theta)+d\}+O(n^{-1})$ and $\kappa_3=-6d+O(n^{-1})$, while for $T(\psi)=M_1(\bar\psi)\{-\bar M_{11}(\bar\psi)\}^{-1/2}$, we have $\kappa_1=-\{g_{\rm INF}(\theta)-2d\}+O(n^{-1})$ and $\kappa_3=12d+O(n^{-1})$. 

Implementation of higher-order inference to error of order $O(n^{-1})$ requires that we estimate the adjustment terms 
${1\over 6}\kappa_3$ and $-\kappa_1+{1 \over 6}\kappa_3$; we might, for example, use plug-in estimates or derive estimates from a simulation procedure such as the parametric bootstrap. If these adjustment terms change rapidly with the value of the 
parameter $\theta$, then there is greater scope for error in the estimation process than if possible if the adjustment terms are stable across $\theta$ values. This observation points to the use of asymptotically standard normal pivots $T(\psi)$ that are derived from the adjusted profile log-likelihood function, since the adjustment terms for such pivots depend only on $g_{\rm INF}(\theta)$ and $d$. If the adjustment terms are small in magnitude, then they are unlikely to vary unduly with $\theta$, and the adjustments can be estimated more reliably. Situations can arise, as is the case in the normal regression example, that the quantity $g_{\rm NP}(\theta)$ is large yet it remains constant with respect to $\theta$. In these circumstances, the need to use the adjusted profile log-likelihood is not so pressing; indeed, for the normal regression model, the parametric bootstrap affords exact inferences, except for simulation error. Since such situations are not commonplace, there is strong motivation for using generally procedures that ensure the magnitudes of the adjustment terms are controlled. However, it could be useful to develop conditions that easily identify models, such as the normal linear regression model, for which the adjustment terms, especially $g_{\rm NP}(\theta)$, are constant or nearly so, since, in such models, the benefit of using adjusted profile likelihood for accurate inference is not so pronounced and procedures based on the regular profile likelihood are likely to suffice.

\section{Discussion}

Accurate inference on a scalar interest parameter $\psi$ in the presence of a nuisance parameter may be obtained using the signed root likelihood ratio statistic $R(\psi)$. A computationally intensive, but analytically simple, approach bases the inference on a bootstrap estimate of the sampling distribution of $R(\psi)$, constructed by fixing the nuisance parameter at its observed constrained maximum likelihood value. Alternatively, inference can be based on a standard normal approximation to the sampling distribution of an analytically adjusted version of $R(\psi)$. For this latter approach, the gold standard is represented by Barndorff-Nielsen's $R^*$ statistic. The adjustment made by this statistic may be decomposed into a sum of two terms.
These adjustments ${\rm INF}(\psi)$ and ${\rm NP}(\psi)$  are determined to second order, $O_p(n^{-1})$, by their expectations.

We have provided an explicit evaluation of these expectations, allowing new theoretical interpretation of the relative importance of the two adjustments and to the intrinsic difficulty of the inference problem within any specified model.

In particular, quantifying the dependence of the expectations on the nuisance parameter provides
insight to circumstances where the bootstrap and analytic approaches might be expected to perform well in terms of accuracy, even in high dimensional problems and with small sample sizes. We have demonstrated that within a particular model, the importance of the nuisance parameter adjustment may depend not only on the structure of the model, as expressed by the nuisance parameter dimension, but the parameter values themselves. In key problems, dependence lies only on the parameter dimension. Calculation of the approximations
$g_{\rm INF}(\theta)$ and $g_{\rm NP}(\theta)$ of $E\{{\rm INF}(\psi)\}$ and $E\{{\rm NP}(\psi)\}$ involves only evaluation of expectations of low order log-likelihood derivatives, and has been demonstrated to give useful theoretical insight to the degree of the adjustment to the signed root statistic $R(\psi)$ given by the statistic $R^*(\psi)$ for any specified inference problem, and therefore to the likely value in use of $R^*(\psi)$ or bootstrapping as a means of improving accuracy.

We note that empirical estimation of the means, through the bootstrap principle of estimation of the nuisance parameter, furnishes a simple procedure for adjustment of the signed root likelihood ratio statistic.  A thorough analysis of this empirical adjustment method for the purposes of inference with higher-order accuracy, as well as a comparison of such an empirical adjustment method with alternative approximations, is beyond the scope of this paper.

\section*{Appendix}

\subsection*{Adjusted profile likelihood}

There have been many suggestions to replace the usual profile likelihood function $M(\psi)$ by an adjusted version $\bar M(\psi)=M(\psi)+B(\psi)$, where $B(\psi)$ is an adjustment function whose derivatives with respect to $\psi$ are of order $O_p(1)$.
The likelihood ratio statistic based on the adjusted profile likelihood is $\bar W(\psi)=2\{\bar M(\bar\psi)-\bar M(\psi)\}$, where $\bar\psi$ is the point at which $\bar M(\psi)$ is maximized.
The signed root of the likelihood ratio statistic based on the adjusted profile likelihood is $\bar R(\psi)={\rm sgn}(\bar\psi-\psi)\{\bar W(\psi)\}^{1/2}$.

Following our previous notation, we write $B_1(\psi)=\partial B(\psi)/\partial\psi$, $B_{11}(\psi)=\partial^2 B(\psi)/\partial\psi^2$, etc.
Let $\beta_1=E\{B_1(\psi)\}$, $\beta_{11}=E(B_{11})$, etc.; these quantities are assumed to be of order $O(1)$.
Further, let $b_1=B_1(\psi)-\beta_1$, $b_{11}=B_{11}(\psi)-\beta_{11}$, etc., with these quantities assumed to be of order $O_p(n^{-1/2})$.
Assume also that the joint cumulants of $nb_1$, $nb_{11}$, $l_r$, $l_{rs}$, etc.\ are of order $O(n)$.

In many instances, the adjustment function $B(\psi)$ has been proposed to take into account the effect of nuisance parameters for inference about $\psi$; see, notably, Cox \& Reid (1987), Barndorff-Nielsen (1983), Skovgaard (1996), Severini (1998), DiCiccio \& Martin (1993), Barndorff-Nielsen \& Chamberlin (1994). These adjustment functions have the effect of reducing the expectation of the profile score from order $O(1)$ to order $O(n^{-1})$.
Specifically, these functions have $\beta_1=\rho+O(n^{-1})$, where $\rho=-\eta\lambda^{1r}\nu^{st}({\textstyle{1\over 2}}\lambda_{rst}+\lambda_{rs,t})$. Since, in general, $E\{M_1(\psi)\}=-\rho+O(n^{-1})$, it follows that $E\{\bar M_1(\psi)\}=O(n^{-1})$: see McCullagh \& Tibshirani (1990), DiCiccio et al.\ (1996).

For a general adjustment function $B(\psi)$, DiCiccio \& Stern (1994b) showed that $\bar R(\psi)=\eta^{1/2}\{\bar R_1+\bar R_2+O_p(n^{-3/2})\}$, where $\bar R_1=R_1=-\lambda^{1r}l_r$ and $\bar R_2=R_2-\lambda^{11}\beta_1$; in particular, $\bar R(\psi)=R(\psi)+\eta^{-1/2}\beta_1+O_p(n^{-1})$. Below, we use this result with a particular adjustment function to obtain a representation of the nuisance parameter adjustment ${\rm NP}(\psi)$, from which $E\{{\rm NP}(\psi)\}$ is then
determined to $O(n^{-1})$. Combined with (\ref{expectR}), this enables calculation to $O(n^{-1})$ of $E\{{\rm INF}(\psi)\}$.

\subsection*{Expectations of Adjustments}

We have,
\begin{eqnarray}
\lefteqn{E\{{\rm NP}(\psi)\}+E\{{\rm INF}(\psi)\}=-E\{R(\psi)\}+O(n^{-1})} \nonumber  \\
&\hspace{-0.5cm} =-\eta^{1/2}\{\lambda^{1r}\lambda^{st}\lambda_{rs,t}
+{\textstyle{1\over 2}}\lambda^{1r}\tau^{st}\lambda_{rs,t}
+{\textstyle{1\over 2}}\lambda^{1r}\lambda^{st}\lambda_{rst}
+{\textstyle{1\over 3}}\lambda^{1r}\tau^{st}\lambda_{rst}\}+O(n^{-1}).
\label{ER}
\end{eqnarray}
It is easily seen that ${\rm NP}(\psi)$ and ${\rm INF}(\psi)$ are of the form ${\rm NP}(\psi)=E\{{\rm NP}(\psi)\}+O_p(n^{-1})$ and ${\rm INF}(\psi)=E\{{\rm INF}(\psi)\}+O_p(n^{-1})$.
Here we develop explicit approximations for $E\{{\rm NP}(\psi)\}$ and $E\{{\rm INF}(\psi)\}$.

The quantity ${\rm NP}(\psi)$ is related to the modified profile likelihood of Barndorff-Nielsen (1983), an adjusted profile likelihood which reduces the bias of the profile score. Following Sartori et al.\ (1999) and Pierce \& Bellio (2006), we have that, up to an additive constant, the log modified profile likelihood is
\begin{align*}
L^{MP}(\psi)&=-R(\psi){\rm NP}(\psi)-\{R(\psi)\}^2/2 \\
&=-R(\psi){\rm NP}(\psi)-M(\hat\psi)+M(\psi) \\
&=-\frac{1}{2}\{R(\psi)+{\rm NP}(\psi)\}^2+O_p(n^{-1}).
 \end{align*}
 The modified profile likelihood therefore corresponds to an adjustment function of the form $B(\psi)=-R(\psi){\rm NP}(\psi)$. Further, the signed square root of the modified profile likelihood ratio statistic is equivalent, to $O_p(n^{-1})$, to $R(\psi)+{\rm NP}(\psi)$, as noted by Sartori et al.\ (1999). The general result of DiCiccio \& Stern (1994b) then gives ${\rm NP}(\psi)=\eta^{-1/2}\beta_1+O_p(n^{-1})$.

Observing that $R(\psi)=(\hat\psi-\psi)\hat\eta^{1/2}+O_p(n^{-1/2})$ and ${\rm NP}(\psi)={\rm NP}(\hat\psi)+O_p(n^{-1})$, we have
\begin{align*}
L^{MP}(\psi)&=-R(\psi){\rm NP}(\psi)-M(\hat\psi)+M(\psi) \\
&=(\psi-\hat\psi)\hat\eta^{1/2}{\rm NP}(\hat\psi)-M(\hat\psi)+M(\psi)+O_p(n^{-1}),
\end{align*}
and differentiation with respect to $\psi$ yields
\[
L^{MP}_1(\psi)=\hat\eta^{1/2}{\rm NP}(\hat\psi)+M_1(\psi)+O_p(n^{-1/2})=\eta^{1/2}{\rm NP}(\psi)+M_1(\psi)+O_p(n^{-1/2}).
\]

Since (see, for example, DiCiccio et al., 1996) $E\{L^{MP}_1(\psi)\}=O(n^{-1})$ and $E\{M_1(\psi)\}=-\rho+O(n^{-1/2})$, it follows that
\[
E\{{\rm NP}(\psi)\}=\eta^{-1/2}\rho+O(n^{-1})=-\eta^{1/2}\lambda^{1r}\nu^{st}(\lambda_{rs,t}+{\textstyle{1\over 2}}\lambda_{rst})+O(n^{-1}),
\]
so that $\beta_1=\eta^{1/2}E\{{\rm NP}(\psi)\}+O(n^{-1/2})=-\eta\lambda^{1r}\nu^{st}(\lambda_{rs,t}+{\textstyle{1\over 2}}\lambda_{rst})+O(n^{-1/2})$.
It follows from (\ref{ER}) that $$E\{{\rm INF}(\psi)\}=\eta^{1/2}\lambda^{1r}\tau^{st}({\textstyle{1\over 2}}\lambda_{rs,t}+{\textstyle{1\over 6}}\lambda_{rst})+O(n^{-1}).$$ We observe also that this analysis confirms $E\{{\rm NP}(\psi)\}=\eta^{-1/2}\beta_1+O(n^{-1})={\rm NP}(\psi)+O_p(n^{-1})$, as noted earlier.

\section*{References}
\begin{enumerate}
\item
{\sc Barndorff-Nielsen, O. E.} (1983)
\newblock On a formula for the conditional distribution of the maximum likelihood estimator.
\newblock {\em Biometrika} {\bf 70}, 343--65.
\item
{\sc Barndorff-Nielsen, O. E.} (1986)
\newblock Inference on full or partial parameters based on the standardized signed log likelihood ratio.
\newblock {\em Biometrika} {\bf 73}, 307--22.
\item
{\sc Barndorff-Nielsen, O. E. \& Chamberlin, S. R.} (1994)
\newblock Stable and invariant adjusted directed likelihoods.
\newblock {\em Biometrika} {\bf 81}, 485--99.
\item
{\sc Barndorff-Nielsen, O. E. \& Cox, D. R.} (1994)
\newblock {\em Inference and Asymptotics}.
\newblock London: Chapman and Hall.
\item
{\sc Brazzale, A. R. \& Davison, A. C.} (2008)
\newblock Accurate parametric inference for small samples. {\it Stat.\ Sci.\ } {\bf 23}, 465--84.
\item
{\sc Brazzale, A. R., Davison, A. C. \& Reid, N.} (2007)
\newblock {\em Applied Asymptotics: Case Studies in Small-Sample Statistics}.
\newblock Cambridge: Cambridge University Press.
\item
{\sc Cox, D. R. \& Reid, N.} (1987)
\newblock Parameter orthogonality and approximate conditional inference (with discussion).
\newblock {\em J.R. Statist. Soc}. B {\bf 53}, 79--109.
\item
{\sc DiCiccio, T. J. \& Efron, B.} (1996)
\newblock Bootstrap confidence intervals (with discussion).
\newblock {\em Stat. Sci. } {\bf 11}, 189--228.
\item
{\sc DiCiccio, T. J. \& Martin, M. A.} (1993)
\newblock Simple modifications for signed roots of likelihood ratio statistics.
\newblock {\em J.R. Statist. Soc}. B {\bf 55}, 305--16.
\item
{\sc DiCiccio, T. J. \& Stern, S. E.} (1994a)
\newblock Frequentist and Bayesian Bartlett correction of test statistics based on adjusted profile likelihoods.
\newblock {\em J.R. Statist. Soc}. B {\bf 56}, 397--408.
\item
{\sc DiCiccio, T. J. \& Stern, S. E.} (1994b)
\newblock Constructing approximately standard normal pivots from signed roots of adjusted likelihood ratio statistics.
\newblock {\em Scand. J. Statist.} {\bf 21}, 447--60.
\item
{\sc DiCiccio, T. J. \& Young, G. A.} (2008)
\newblock Conditional properties of unconditional parametric bootstrap procedures for inference in exponential families.
\newblock {\em Biometrika} {\bf 95}, 747--58.
\item
{\sc DiCiccio, T. J., Martin, M. A. \& Stern, S. E.} (2001)
\newblock Simple and accurate one-sided inference from signed roots of likelihood ratios.
\newblock {\em Can. J. Statist.} {\bf 29}, 67--76.
\item
{\sc DiCiccio, T. J., Martin, M. A., Stern, S. E. \& Young, G. A.} (1996)
\newblock Information bias and adjusted
profile likelihood.
\newblock {\em J.R. Statist. Soc}. B {\bf 58}, 189--203.
\item
{\sc DiCiccio, T. J., Kuffner, T. A., Young, G. A. \& Zaretzki, R.} (2015)
\newblock Stability and uniqueness of p-values for likelihood-based inference.
\newblock To appear in  {\em Statistica Sinica}.
\item
{\sc Efron, B.} (1987)
\newblock Better bootstrap confidence intervals (with discussion).
\newblock {\em J. Amer. Statist. Assoc.} {\bf 82}, 171--200.
\item
{\sc Fraser, D. A. S.} (1990)
\newblock Tail probabilities from observed likelihoods.
\newblock {\em Biometrika} {\bf 77}, 65--76.
\item
{\sc Lawley, D. N.} (1956)
\newblock A general method for approximating to the distribution of likelihood ratio criteria.
\newblock {\em Biometrika} {\bf 43}, 295--303.
\item
{\sc Lee, S. M. S. \& Young, G. A.} (2005)
\newblock Parametric bootstrapping with nuisance parameters.
\newblock {\em Stat. Prob. Letters} {\bf 71}, 143--53.
\item
{\sc McCullagh, P. \& Tibshirani, R.} (1990)
\newblock A simple method for the adjustment of profile likelihoods.
 \newblock {\em J. Roy. Statist. Soc. B} {\bf 52}, 325--44.
\item
{\sc Pierce, D. A. \& Bellio, R.} (2006)
\newblock Effects of the reference set on frequentist inferences.
\newblock {\em Biometrika} {\bf 93}, 425�-38.
\item
{\sc Pierce, D. A. \& Peters, D.} (1992)
\newblock Practical use of higher-order asymptotics for multiparameter
exponential families (with discussion).
\newblock {\em J.R. Statist. Soc}. B {\bf 54}, 701--38.
\item
{\sc Sartori, N.} (2003)
\newblock Modified profile likelihoods in models with stratum nuisance parameters.
\newblock {\em Biometrika} {\bf 90}, 533--49.
\item
{\sc Sartori, N., Bellio, R., Salvan, A. \& Pace, L.} (1999)
\newblock The directed modified profile likelihood in models with many nuisance parameters.
\newblock {\em Biometrika} {\bf 86}, 735--42.
\item
{\sc Severini, T. A.} (1998)
\newblock An approximation to the modified profile likelihood function.
\newblock {\em Biometrika} {\bf 85}, 403--11.
\item
{\sc Severini, T. A.} (2000)
\newblock {\em Likelihood methods in Statistics}.
\newblock Oxford: Oxford University Press.
\item
{\sc Skovgaard, I. M.} (1996)
\newblock An explicit large-deviation approximation to one-parameter tests.
\newblock {\em Bernoulli} {\bf 2}, 145--65.
\item
{\sc Young, G. A.} (2009)
\newblock Routes to higher-order accuracy in parametric inference.
\newblock {\em Aust. N.Z. J. Stat.} {\bf 51}, 115--26.
\end{enumerate}

\end{document}